\documentclass[a4paper,12pt]{amsart}
\usepackage[latin1]{inputenc}
\usepackage[T1]{fontenc}
\usepackage{amsfonts,amssymb}
\usepackage{pstricks}
\usepackage{amsthm}
\usepackage[all]{xy}
\usepackage[pdftex]{graphicx}
\usepackage{graphicx}  
 \usepackage[first=0,last=9]{lcg}
 \usepackage{tikz}
 \usetikzlibrary{tikzmark}
\usetikzlibrary{calc, positioning, fit, shapes.misc}

 \usepackage{color, colortbl}

 \tikzset{every picture/.style={remember picture}}

 \definecolor{Gray}{gray}{0.9}
\newtheorem*{Theo}{Theorem}

\newtheorem{theo}{Theorem}[section]
\newtheorem{defi}{Definition}[section]
\newtheorem{prop}[theo]{Proposition}

\newtheorem{lem}[theo]{Lemma}

\newtheorem{rem}[theo]{Remark}
\newtheorem{ex}{Example}

\title[Neighborly partitions and  Gordon's identities]{Neighborly partitions, Hypergraphs\\ and Gordon's identities}
\author{Pooneh Afsharijoo and Hussein Mourtada}

\begin {document}

\maketitle

\begin{abstract} We prove a family of partition identities which is "dual" to the family of Andrews-Gordon's identities. These identities are inspired by a correspondence between   a special type of partitions and "hypergraphs" and their proof uses combinatorial commutative algebra.
\end{abstract}

\footnote{keywords: Andrew-Gordon's identities; integer partitions; graded algebras; Hypergraphs; Hilbert series.}


\section{Introduction}
The Andrews-Gordon identities (Andrews 1974 \cite{AG}) are the $q-$series identities which state that for the integers $r$ and $i$ satisfying $ r\geq 2, \ 1 \leq i \leq r,$ we have

\begin{equation}\label{GA}
\sum_{n_1,n_2,\dots n_{r-1} \geq 0} \frac{q^{N_1^2+N_2^2+\dots +N_{r-1}^2+N_i+ N_{i+1}+\dots+N_{r-1}}}{(q)_{n_1}(q)_{n_2}\dots (q)_{n_{r-1}}}=\frac{\prod\limits_{n\geq1,~~ n \equiv 0, \pm i (mod. 2r+1)}(1-q^n)}{\prod\limits_{n\geq1}(1-q^n)}.
\end{equation}
Where $q$ is a variable and $N_j=n_j +n_{j+1} +\dots + n_{r-1}$ for all $1 \leq j\leq r-1 $ and $(q)_n=(1-q)(1-q^2)\cdots(1-q^n)$. In the literature, the right member of the identity (\ref{GA}) is written with the obvious simplification made; we write it in this way to emphasize on the numerator which plays an important role in our paper.  
The Andrews-Gordon identities are generalizations of the famous Rogers-Ramanujan identities that we obtain if we put  $r=2$ and $i=1,2.$ \\

There is a combinatorial (versus analytic) version of the identity $(\ref{GA})$ which is stated in terms of integer partition. Recall that an integer \textit{partition} (of length $\ell$) of a positive integer $n$ is a sequence $\lambda=(\lambda_1\geq \cdots \geq \lambda_\ell)$ of positive integers $\lambda_i$, for $1 \leq i \leq \ell$, such that
$$\lambda_1+\cdots +\lambda_\ell=n.$$
The integers $\lambda_i$ are called the \textit{parts of the partition $\lambda$}.
\\
The combinatorial version of the identities (\ref{GA}) states the following (see Theorem $1$ in \cite{G}):

\begin{Theo} \label{Gordon}  (\textit{Gordon's identities}). Given integers $r\geq 2$ and $1\leq i \leq r,$ let $B_{r,i}(n)$ denote the number of partitions of $n$ of the form $(b_1,\dots, b_s)$, where $b_{j}-b_{j+r-1} \geq 2$ and at most $i-1$ of the integers $b_j$ are equal to $1$. Let $A_{r,i}(n)$ denote the number of partitions of $n$ into parts $\not\equiv 0,\pm i \ ( \text{mod}.  2r+1)$. Then $A_{r,i}(n)=B_{r,i}(n)$ for all integers $n$.

\end{Theo}

Beside combinatorics and number theory, these identities appeared also in representation theory \cite{LZ} and in Algebraic Geometry and Commutative Algebra 
\cite{BMS,Af,ADJM1,ADJM2,M_Hand}).\\

In this paper, we will prove a family of identities which is in some sense (that will be clarified in a moment) dual to Gordon's identities. For that we will introduce the notions of Neighborly partitions and their signature. These notions generalizes the notions with the same name  which were introduced in \cite{MM} to prove dual identities to those of Rogers-Ramanujan. \\

We begin by introducing neighborly partitions. Recall that for an integer partition $\lambda,$ the \textit{multiplicity} $m_\lambda(\lambda_j)$ of a part $\lambda_j$ of $\lambda$ is the number of occurrences of $\lambda_j$ in $\lambda$; for example if $\lambda$ is the partition $3+2+2+1+1+1$ of $10,$ then the multiplicities of $3,2$ and $1$ are respectively  $1,2$ and $3.$

\begin{defi} (Neighborly partitions) Let $\lambda=(\lambda_1\geq \cdots \geq \lambda_m)$ be an integer partition and let $i$ and $r$ be partitions such that $1\leq i \leq r$. We say that $\lambda$ is an $(r,i)$-Neighborly partition if it satisfies the following conditions:
\begin{itemize}
\item[1.] For each part $\lambda_j \neq 1$ we have $1 \leq m_{\lambda}(\lambda_j)\leq r$ and $1 \leq m_{\lambda}(1)\leq i.$
\\
\item[2.] 

If $m_{\lambda}(1)=i,$ then for all parts $\lambda_j \neq 1$ of $\lambda$ there exists a sub-partition $B_j=(\lambda_k\geq \cdots \geq \lambda_{k+r-1})$ of length $r$ of $\lambda$ containing $\lambda_j$ in which $\lambda_k-\lambda_{k+r-1}\leq 1.$
\\
\item[3.] If $m_{\lambda}(1)<i,$ then for all parts $\lambda_j $ of $\lambda$ there exists a sub-partition $B_j=(\lambda_k\geq \cdots \geq \lambda_{k+r-1})$ of length $r$ of $\lambda$ containing $\lambda_j$ in which $\lambda_k-\lambda_{k+r-1}\leq 1.$

\end{itemize}

\end{defi}

We denote the set of $(r,i)$-Neighborly partitions by $\mathcal{N}_{r,i}.$

\begin{rem}\label{B_k}
Note that each sub-partition of type $B_j$ of a $(r,i)$-neighborly partition $\lambda$ is of the form: $$(\underbrace{(\ell+1),\cdots,(\ell+1)}_{\text {(r-s) times}},\underbrace{\ell,\cdots,\ell}_{\text {s-times}}),$$ for some $1\leq s \leq r$ and $\ell \geq 1.$ In particular,  $\lambda_j \in \{l,l+1\},$ has many neighbors; these are the parts $\lambda_k$ satisfying 
$\mid \lambda_k-\lambda_j\mid \leq 1.$ This explains the name neighborly. One remarks that defining conditions of neighborly partitions are of opposite nature to the defining conditions of the partitions $B_{r,i}$ that appear in Gordon's identities. 
\end{rem}



\begin{ex}\label{partition} The integer $5$ has the following partitions:

$$5=4+1=3+2=3+1+1=2+2+1=2+1+1+1=1+1+1+1+1.$$

For $r=3$ we have:

$$\mathcal{N}_{3,3}(5)=\{2+2+1, 2+1+1+1\},\  \mathcal{N}_{3,2}(5)=\{2+2+1\},\ \mathcal{N}_{3,1}(5)=\emptyset.$$
\end{ex}

To a neighborly partition $\lambda,$ we will associate a hypergraph $\mathcal{H}_\lambda$ and a signature which is a number that we define in terms of $\mathcal{H}_\lambda.$\\
Recall that the notion "hypergrpah" is a generalization of the notion "graph" where edges may join more than two vertices. More precisely, a \textit{hypergraph} $\mathcal{H}$ is a pair $(V(\mathcal{H}), E(\mathcal{H}))$ where $V(\mathcal{H})$ is a set of elements called the vertices of $\mathcal{H}$ and $E(\mathcal{H})$ is a set of subsets of $V(\mathcal{H})$ called the edges of $\mathcal{H}$.\\
One can represent a hypergraph $\mathcal{H}$ graphically by a \textit{Parallel Aggregated Ordered Hypergraph (PAOH for short)}; A $\textit{PAOH}$ represents the vertices of $\mathcal{H}$ by parallel horizontal rows, and the edges of $\mathcal{H}$ by vertical lines in which a point represents a vertex of the edge (see Example \ref{first}). A hypergraph is called \textit{$k$-uniform} if each edge contains exactly $k$ vertices. Thus, a $2$-uniform hypergraph is a graph. \\
The \textit{degree}, $\deg(v),$ of a vertex $v$ in a hypergraph $\mathcal{H}$ is the number of edges of $\mathcal{H}$ containing $v.$ If $\deg(v)=0$ then we say that $v$ is an \textit{isolated vertices} of $\mathcal{H}.$ 
A hypergraph $\mathcal{H}$ is simple if there is no edge of $\mathcal{H}$ which contains another edge.\\
A \textit{vertex induced sub-hypergraph} $\mathcal{L}$ of $\mathcal{H}$ is a hypergraph whose set of vertices $V(\mathcal{L})$ is a subset of $V(\mathcal{H})$ and its edges are the edges of $\mathcal{H}$ whose vertices are in $V(\mathcal{L}).$\\
An \textit{edge induced sub-hypergraph} $\mathcal{L}$ of $\mathcal{H}$ is a hypergraph whose edge set is a sub set of $E(\mathcal{H})$ and whose vertex set is the union of the vertices of its edges.

\begin{ex}\label{first} Consider a hypergraph $\mathcal{H}$ with $V(\mathcal{H})=\{v_1,\cdots,v_6\}$ and  $E(\mathcal{H})=\{(v_1,v_2,v_3,v_4),(v_1,v_3),(v_2,v_5)\}.$ Then its \textit{PAOH} is represented as follows:

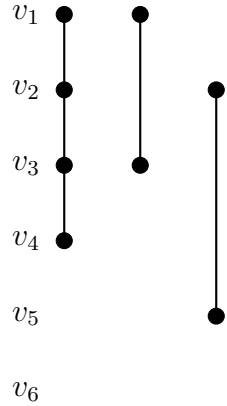
\begin{center}

\begin{tikzpicture}
\draw[fill=black] (0,0) circle (3pt);
\draw[fill=black] (0,-1) circle (3pt);
\draw[fill=black] (0,-2) circle (3pt);
\draw[fill=black] (0,-3) circle (3pt);
\draw[fill=black] (1,0) circle (3pt);
\draw[fill=black] (1,-2) circle (3pt);
\draw[fill=black] (2,-1) circle (3pt);
\draw[fill=black] (2,-4) circle (3pt);

\node at (-0.5,0) {$v_1$};
\node at (-0.5,-1) {$v_2$};
\node at (-0.5,-2) {$v_3$};
\node at (-0.5,-3) {$v_4$};
\node at (-0.5,-4) {$v_5$};
\node at (-0.5,-5) {$v_6$};
\draw[thick] (0,0) -- (0,-1);
\draw[thick] (0,-1) -- (0,-2);
\draw[thick] (0,-2) -- (0,-3);
\draw[thick] (1,0) -- (1,-2);
\draw[thick] (2,-1) -- (2,-4);

\end{tikzpicture}
\end{center}
and we have:
\begin{itemize}
\item[-]$\deg(v_1)=\deg(v_2)=\deg(v_3)=2, \ \deg(v_4)=\deg(v_5)=1$ and $v_6$ is an isolated vertices of $\mathcal{H}.$
\\
\item[-] The hypergraph $\mathcal{H}$ is not simple since $(v_1,v_3)\subset (v_1,v_2,v_3,v_4).$ 
\\
\item[-] The hypergraph $\mathcal{L}_1$ with $V(\mathcal{L}_1)=\{v_1,v_2,v_3,v_5,v_6\}$ and $E(\mathcal{L}_1)=\{(v_1,v_3),(v_2,v_5)\}$ is a vertex induced sub-hypergraph of $\mathcal{H}.$
\\
\item[-] The hypergraph $\mathcal{L}_2$ with  $E(\mathcal{L}_2)=\{(v_1,v_3),(v_2,v_5)\}$ and $V(\mathcal{L}_2)=\{v_1,v_2,v_3,v_5\}$ is an edge induced sub-hypergraph of $\mathcal{H}.$ 
\end{itemize}
\end{ex}

To each $(r,i)-$neighborly partition $\lambda$ we associate an hypergraph $\mathcal{H}_\lambda$ as follows: 

\noindent The set of vertices of $\mathcal{H}_\lambda$ is in bijection with the parts of $\lambda$:
$$V(\mathcal{H}_\lambda)=\{x_{j,k}| \ j \ \text{is a part of} \ \lambda \ \text{and} \ 1\leq k \leq m_{\lambda}(j) \}.$$
\noindent The set edges of $\mathcal{H}_\lambda$ is in  bijection with the set of all sub-partitions of type $B_j$ of $\lambda$ (see Remark \ref{B_k}), so we have:

\begin{itemize}

\item if $m_{\lambda}(1)<i,E(\mathcal{H}_\lambda)=\{ (x_{\ell,1},\cdots,x_{\ell,s},x_{(\ell+1),1},\cdots,x_{(\ell+1),(r-s)})| \ \text{for all} \ x_{j,k}\in V(\mathcal{H}_\lambda) \ \text{and} \ 1\leq s \leq r\}.$


\item If $m_{\lambda}(1)=i$ then we add to the above set of edges the edge $(x_{1,1},\cdots,x_{1,i}).$

\end{itemize}

Note that if $\ell$ is a part of $\lambda$ with $m(\ell)=r$ then the edge associated to the sub-partition  $(\underbrace{\ell,\cdots,\ell}_{\text {r-times}})$  of $\lambda$ is $(x_{\ell,1},\cdots,x_{\ell,r})$. This is the case $s=r.$


Note also that if $i=r$ or if $1\leq i <r$ and $0 \leq m_{\lambda}(1)<i$ then $\mathcal{H}_\lambda$ is a $r-$uniform hypergraph.

\begin{ex} To the partition $\lambda= 2+1+1+1$ of $\mathcal{N}_{3,3}(5)$ we associate a hypergraph $\mathcal{H}_{\lambda}$ whose vertex set and edge set are:
$$V(\mathcal{H}_\lambda)=\{x_{1,1}, x_{1,2}, x_{1,3}, x_{2,1} \}, E(\mathcal{H}_\lambda)=\{(x_{1,1},x_{1,2},x_{1,3}),(x_{1,1}, x_{1,2}, x_{2,1})\}.$$

The PAOH representation of $\lambda$ is as follows:
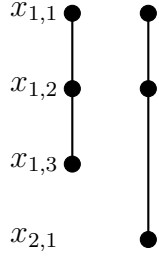
\begin{center}

\begin{tikzpicture}
\draw[fill=black] (0,0) circle (3pt);
\draw[fill=black] (0,-1) circle (3pt);
\draw[fill=black] (0,-2) circle (3pt);
\draw[fill=black] (1,0) circle (3pt);
\draw[fill=black] (1,-1) circle (3pt);
\draw[fill=black] (1,-3) circle (3pt);

\node at (-0.5,0) {$x_{1,1}$};
\node at (-0.5,-1) {$x_{1,2}$};
\node at (-0.5,-2) {$x_{1,3}$};
\node at (-0.5,-3) {$x_{2,1}$};
\draw[thick] (0,0) -- (0,-1);
\draw[thick] (0,-1) -- (0,-2);
\draw[thick] (1,0) -- (1,-1);
\draw[thick] (1,-1) -- (1,-3);
\end{tikzpicture}
\end{center}
\end{ex}





   

For a hypergraph $\mathcal{H},$ we denote by $Sub_v(\mathcal{H})$ (respectively by $Sub_e(\mathcal{H})$) the set of \textbf{vertex induced sub-hypergraphs} (respectively \textbf{edge induced sub-hypergraphs}) of $\mathcal{H}$ \textbf{without isolated vertices}. We will also denote by $|E(\mathcal{H})|$ the number of edges of  $\mathcal{H}$.

\begin{defi} Let $\lambda \in \mathcal{N}_{r,i}.$ We define the signature of $\lambda$ as follows:
$$\delta(\lambda)=\sum_{\substack{\mathcal{L} \in Sub_e(\mathcal{H}_{\lambda}) \\ V(\mathcal{L})=V(\mathcal{H}_{\lambda})}} (-1)^{|E(\mathcal{L})|},$$
\end{defi}

Let us now to denote by $\mathcal{R}_{r,i}(n)$ the set of integer partitions of $n$ whose parts are distinct, congruent to $0,-i$ or $i$ modulo $2r+1$. The main result of this paper is the following theorem (see Theorem \ref{generatrice}):

\begin{theo}\label{product} Let $1\leq i \leq r$ be the integers. Then:

$$\sum_{\lambda \in \mathcal{N}_{r,i}}\delta(\lambda) q^{|\lambda|}= \sum_{n\in \mathbb{N}}\mathcal{R}_{r,i}(n)q^n= \prod_{j \equiv 0, \pm i [2r+1]} (1-q^j),$$
where $|\lambda|$ is the sum of the parts of $\lambda$.
\end{theo}
This is equivalent to the following theorem (see Theorem \ref{coeff.}):

\begin{theo} Let $1\leq i \leq r$ be the integers. Then:
$$\sum_{\lambda \in \mathcal{N}_{r,i}(n)} \delta(\lambda)=\sum_{\lambda\in \mathcal{R}_{r,i}(n)} (-1)^{\ell(\lambda)}.$$
\end{theo}

\begin{ex} For the partitions of the integer $7$ we have:
$$\mathcal{N}_{3,3}(7)=\{\underbrace{3+2+2}_{\alpha}, \underbrace{2+2+2+1}_{\beta}, \underbrace{2+2+1+1+1}_{\gamma}\},\  \mathcal{R}_{3,3}(7)=\{7, 4+3\}.$$
we have:
$$\sum_{\lambda\in \mathcal{R}_{3,3}(7)} (-1)^{\ell(\lambda)}=(-1)^{\ell(7)}+(-1)^{\ell(4+3)}=(-1)^1+(-1)^2=0.$$

The hyper graphs associated to the partitions $\alpha, \beta$ and $\gamma$ are , from left to right, as follows:

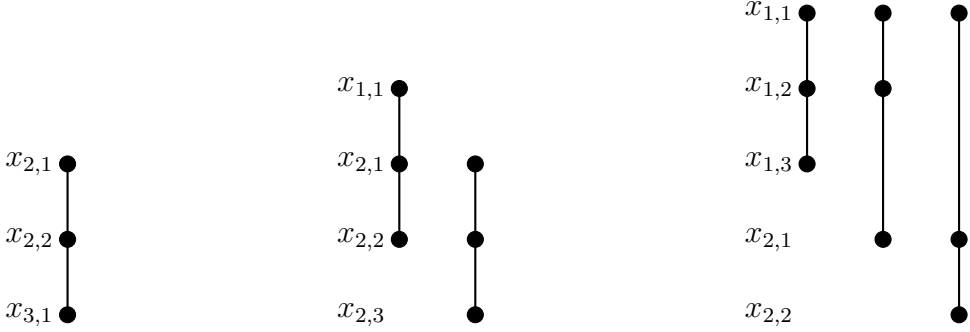
\begin{center}

\begin{tikzpicture}
\draw[fill=black] (0,2) circle (3pt);
\draw[fill=black] (0,1) circle (3pt);
\draw[fill=black] (0,0) circle (3pt);


\node at (-0.5,2) {$x_{2,1}$};
\node at (-0.5,1) {$x_{2,2}$};
\node at (-0.5,0) {$x_{3,1}$};
\draw[thick] (0,2) -- (0,1);
\draw[thick] (0,1) -- (0,0);
\end{tikzpicture}
\hspace{3cm}
\begin{tikzpicture}
\draw[fill=black] (0,0) circle (3pt);
\draw[fill=black] (0,-1) circle (3pt);
\draw[fill=black] (0,-2) circle (3pt);
\draw[fill=black] (1,-1) circle (3pt);
\draw[fill=black] (1,-2) circle (3pt);
\draw[fill=black] (1,-3) circle (3pt);


\node at (-0.5,0) {$x_{1,1}$};
\node at (-0.5,-1) {$x_{2,1}$};
\node at (-0.5,-2) {$x_{2,2}$};
\node at (-0.5,-3) {$x_{2,3}$};

\draw[thick] (0,0) -- (0,-1);
\draw[thick] (0,-1) -- (0,-2);
\draw[thick] (1,-1) -- (1,-2);
\draw[thick] (1,-2) -- (1,-3);
\end{tikzpicture}
\hspace{3cm}
\begin{tikzpicture}
\draw[fill=black] (0,0) circle (3pt);
\draw[fill=black] (0,-1) circle (3pt);
\draw[fill=black] (0,-2) circle (3pt);
\draw[fill=black] (1,0) circle (3pt);
\draw[fill=black] (1,-1) circle (3pt);
\draw[fill=black] (1,-3) circle (3pt);
\draw[fill=black] (2,0) circle (3pt);
\draw[fill=black] (2,-3) circle (3pt);
\draw[fill=black] (2,-4) circle (3pt);


\node at (-0.5,0) {$x_{1,1}$};
\node at (-0.5,-1) {$x_{1,2}$};
\node at (-0.5,-2) {$x_{1,3}$};
\node at (-0.5,-3) {$x_{2,1}$};
\node at (-0.5,-4) {$x_{2,2}$};
\draw[thick] (0,0) -- (0,-1);
\draw[thick] (0,-1) -- (0,-2);
\draw[thick] (1,0) -- (1,-1);
\draw[thick] (1,-1) -- (1,-3);
\draw[thick] (2,0) -- (2,-3);
\draw[thick] (2,-3) -- (2,-4);

\end{tikzpicture}
\end{center}

The unique edge induced hypergraph of $\mathcal{H}_{\alpha}$ with the same vertex set as $\mathcal{H}_{\alpha}$ and without isolated vertices is itself. Since it has just one edge, hence $\delta(\alpha)=(-1)^1=-1.$
\\
Similarly, the unique hypergraph in $Sub_e(\mathcal{H}_{\beta})$ with the vertex set $V(\mathcal{\beta})$ is itself. Since it has two edges, we have $\delta(\beta)=(-1)^2=1.$
\\
But, the set $Sub_e(\mathcal{H}_{\gamma})$ contains two hypergraphs with the vertex set equal to $V(\mathcal{H}_{\gamma}):$ the hypergraph $\gamma$ and the following hypergraph which has two edges:
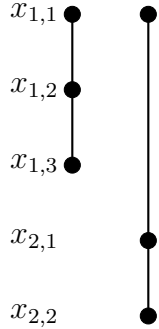
\begin{center}
\begin{tikzpicture}
\draw[fill=black] (0,0) circle (3pt);
\draw[fill=black] (0,-1) circle (3pt);
\draw[fill=black] (0,-2) circle (3pt);
\draw[fill=black] (1,0) circle (3pt);
\draw[fill=black] (1,-3) circle (3pt);
\draw[fill=black] (1,-4) circle (3pt);


\node at (-0.5,0) {$x_{1,1}$};
\node at (-0.5,-1) {$x_{1,2}$};
\node at (-0.5,-2) {$x_{1,3}$};
\node at (-0.5,-3) {$x_{2,1}$};
\node at (-0.5,-4) {$x_{2,2}$};
\draw[thick] (0,0) -- (0,-1);
\draw[thick] (0,-1) -- (0,-2);
\draw[thick] (1,0) -- (1,-3);
\draw[thick] (1,-3) -- (1,-4);
\end{tikzpicture}
\end{center}

Thus, $\delta(\gamma)=(-1)^3+(-1)^2=0$ and therefore:
$$\sum_{\lambda \in \mathcal{N}_{r,i}(7)} \delta(\lambda)=\delta(\alpha)+\delta(\beta)+\delta(\gamma)=-1+1+0=0,$$
which is equal to $\sum_{\lambda\in \mathcal{R}_{r,i}(n)} (-1)^{\ell(\lambda)}$ and the theorem holds.

\end{ex}
The main theorem greatly generalizes the main results of \cite{MM}. It is worth noticing
that the proof not only uses the Andrews-Gordon's identities, but are actually equivalent 
to them, which means that a direct proof of our theorem gives also another proof of the Andrews-Gordon's identities. This program was very recently pursued in the case of Rogers-Ramanujan's identities by O'Hara and Stanton \cite{OS}.  \\

The proof of our main results follows the organization of the paper: in the second section, we introduce an 
infinite hypergraph $\mathcal{H}_{r,i}^{\infty}$ and we express the left member of the identity in theorem 
\ref{product} via a counting series $S(\textbf{v},y)$ of some finite sub-hypergarph of $\mathcal{H}_{r,i}^{\infty}.$ 
In section three, with a simple hypergraph $\mathcal{H},$ we associate a monomial ideal $\mathcal{I}(\mathcal{H})$ 
in a weighted polynomial ring $A$ whose variables are in bijection with the vertices of $\mathcal{H};$ we then 
consider a kind of Hilbert series $H_{\mathcal{H}}$ of this ideal, this is the generating series of the monomials 
in the quotient of $A$ by $ \mathcal{I}(\mathcal{H});$ this latter series is expressed in the same section via $S(\textbf{v},y).$ We consider a specialization of $H_{\mathcal{H}}$ which we link in section four to the left member of the identity in theorem \ref{product}, and in section five to the right member of the same identity, using
Gordon's identities.

\section{$(r,i)$-Neighborly partitions and hypergraphs}
In this section, we introduce an 
infinite hypergraph $\mathcal{H}_{r,i}^{\infty}$ and we  express the left member of the identity in theorem 
\ref{product} via a counting series $S(\textbf{v},y)$ of some finite sub-hypergarph of $\mathcal{H}_{r,i}^{\infty},$
see Lemma \ref{S-series}.\\ 

For $i,r \in \mathbb{N},1\leq i \leq r, $ consider the infinite hypergraph 
$\mathcal{H}_{r,i}^{\infty}$ with:

\begin{itemize}
\item[-] $V(\mathcal{H}_{r,i}^{\infty})=\{x_{1,1},\cdots, x_{1,i},x_{j,k}|\  k\in [|1,r|],  j\in\mathbb{N}^*\setminus \{1\}  \} $.
\\
\item[-] $E(\mathcal{H}_{r,i}^{\infty})=\{ (x_{1,1},\cdots,x_{1,i}),(x_{\ell,1},\cdots,x_{\ell,s},x_{(\ell+1),1},\cdots,x_{(\ell+1),(r-s)})  \}$,
\\
\\
where $\ell\in \mathbb{N}^*.$ If $\ell=1$ then $1\leq s \leq i-1,$ otherwise $1\leq s \leq r.$
\end{itemize} 

The \textit{PAOH} representation of $\mathcal{H}_{r,r}^{\infty}$ is given in figure 1.

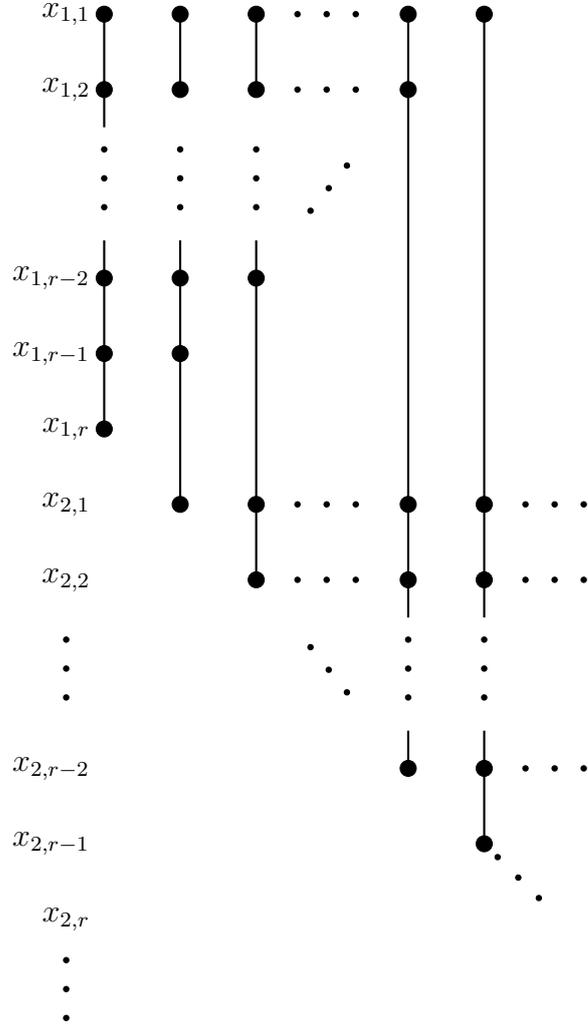
\begin{figure}\caption{The PAOH representation of $\mathcal{H}_{r,r}^{\infty}$}

\begin{center}

\begin{tikzpicture}
\draw[fill=black] (0,0) circle (3pt);
\draw[fill=black] (0,-1) circle (3pt);
\draw[fill=black] (0,-3.5) circle (3pt);
\draw[fill=black] (0,-4.5) circle (3pt);
\draw[fill=black] (0,-5.5) circle (3pt);
\draw[fill=black] (1,0) circle (3pt);
\draw[fill=black] (1,-1) circle (3pt);
\draw[fill=black] (1,-3.5) circle (3pt);
\draw[fill=black] (1,-4.5) circle (3pt);
\draw[fill=black] (1,-6.5) circle (3pt);
\draw[fill=black] (2,0) circle (3pt);
\draw[fill=black] (2,-1) circle (3pt);
\draw[fill=black] (2,-6.5) circle (3pt);
\draw[fill=black] (2,-7.5) circle (3pt);
\draw[fill=black] (2,-3.5) circle (3pt);
\draw[fill=black] (4,0) circle (3pt);
\draw[fill=black] (4,-1) circle (3pt);
\draw[fill=black] (4,-6.5) circle (3pt);
\draw[fill=black] (4,-7.5) circle (3pt);
\draw[fill=black] (4,-10) circle (3pt);
\draw[fill=black] (5,0) circle (3pt);
\draw[fill=black] (5,-6.5) circle (3pt);
\draw[fill=black] (5,-7.5) circle (3pt);
\draw[fill=black] (5,-10) circle (3pt);
\draw[fill=black] (5,-11) circle (3pt);



\node at (-0.5,0) {$x_{1,1}$};
\node at (-0.5,-1) {$x_{1,2}$};
\node at (-0.7,-3.5) {$x_{1,r-2}$};
\node at (-0.7,-4.5) {$x_{1,r-1}$};
\node at (-0.5,-5.5) {$x_{1,r}$};
\node at (-0.5,-6.5) {$x_{2,1}$};
\node at (-0.5,-7.5) {$x_{2,2}$};
\node at (-0.7,-10) {$x_{2,r-2}$};
\node at (-0.7,-11) {$x_{2,r-1}$};
\node at (-0.5,-12) {$x_{2,r}$};

\draw[thick] (0,0) -- (0,-1);
\draw[thick] (0,-1) -- (0,-1.5);
\path (0,-1.5) -- (0,-3) node [black, font=\Huge, midway, sloped] {$\dots$};
\draw[thick] (0,-3) -- (0,-3.5);
\draw[thick] (0,-3.5) -- (0,-4.5);
\draw[thick] (0,-4.5) -- (0,-5.5);
\draw[thick] (1,0) -- (1,-1);
\path (1,-1.5) -- (1,-3) node [black, font=\Huge, midway, sloped] {$\dots$};
\draw[thick] (1,-3) -- (1,-3.5);
\draw[thick] (1,-3.5) -- (1,-4.5);
\draw[thick] (1,-4.5) -- (1,-6.5);
\draw[thick] (2,0) -- (2,-1);
\path (2,-1.5) -- (2,-3) node [black, font=\Huge, midway, sloped] {$\dots$};
\draw[thick] (2,-3) -- (2,-3.5);
\draw[thick] (2,-3.5) -- (2,-6.5);
\draw[thick] (2,-6.5) -- (2,-7.5);
\draw[thick] (4,0) -- (4,-1);
\draw[thick] (4,-1) -- (4,-6.5);
\draw[thick] (4,-6.5) -- (4,-7.5);
\draw[thick] (4,-7.5) -- (4,-8);
\path (4,-8) -- (4,-9.5) node [black, font=\Huge, midway, sloped] {$\dots$};
\draw[thick] (4,-9.5) -- (4,-10);
\draw[thick] (5,0) -- (5,-6.5);
\draw[thick] (5,-6.5) -- (5,-7.5);
\draw[thick] (5,-7.5) -- (5,-8);
\path (5,-8) -- (5,-9.5) node [black, font=\Huge, midway, sloped] {$\dots$};
\draw[thick] (5,-9.5) -- (5,-10);
\draw[thick] (5,-10) -- (5,-11);

\path (2,0) -- (4,0) node [black, font=\Huge, midway, sloped] {$\dots$};

\path (2,-1) -- (4,-1) node [black, font=\Huge, midway, sloped] {$\dots$};
\path (2,-6.5) -- (4,-6.5) node [black, font=\Huge, midway, sloped] {$\dots$};
\path (2,-7.5) -- (4,-7.5) node [black, font=\Huge, midway, sloped] {$\dots$};
\path (5,-6.5) -- (7,-6.5) node [black, font=\Huge, midway, sloped] {$\dots$};
\path (5,-7.5) -- (7,-7.5) node [black, font=\Huge, midway, sloped] {$\dots$};
\path (5,-10) -- (7,-10) node [black, font=\Huge, midway, sloped] {$\dots$};
\path (5,-11) -- (6,-12) node [black, font=\Huge, midway, sloped] {$\dots$};

\path (4,-1) -- (2,-3.5) node [black, font=\Huge, midway, sloped] {$\dots$};
\path (2,-7.5) -- (4,-10) node [black, font=\Huge, midway, sloped] {$\dots$};


\path (-0.5,-7.5) -- (-0.5,-10) node [black, font=\Huge, midway, sloped] {$\dots$};
\path (-0.5,-12.5) -- (-0.5,-13.5) node [black, font=\Huge, midway, sloped] {$\dots$};

\end{tikzpicture}
\end{center}

\end{figure}

Note that, the set of \textbf{finite} sub-hypergrapghes in $Sub_{\textbf{v}}(\mathcal{H}_{r,i}^{\infty})$ is in bijection with the set of integer partitions $\lambda \in \mathcal{N}_{r,i}.$

\begin{defi} Let $\mathcal{H}$ be a simple hypergraph with a countable vertex set $V(\mathcal{H})=\{v_h| \ h\in I\}.$ We define the following multivariable series in $\textbf{v}:=(v_h)_{h \in I}.$ :

$$S(\textbf{v},y)=\sum_{\substack{\mathcal{L} \in Sub_e(\mathcal{H}) \\\text{finite}}} (\prod_{v_h \in V(\mathcal{L})}v_h) \ y^{|E(\mathcal{L})|},$$
 
where we assume that the coefficient of $y^0$ is equal to $1$
\end{defi}
We have the following lemma which gives us an expression of the left hand side series in the Theorem \ref{generatrice} in terms of the series just defined above.

\begin{lem}\label{S-series} Let $\textbf{v}=(x_{j,k})_{x_{j,k} \in V(\mathcal{H}_{r,i}^{\infty})}$ and denote by $S_{r,i}^{w}(q,y)$ the series obtained from $S(\textbf{v},y)$ by replacing each $x_{j,k}$ by $q^j$. Then we have:

$$S_{r,i}^{w}(q,-1)=\sum_{\lambda \in \mathcal{N}_{r,i}} \delta(\lambda) q^{|\lambda|}.$$

\end{lem} 

\begin{proof} By definition of $S_\mathcal{H}(\textbf{v},y)$ we have:
\begin{align*}
 S(\textbf{v},y)=& \sum_{\substack{\mathcal{L} \in Sub_e(\mathcal{H}_{r,i}^{\infty}) \\\text{finite}}} (\prod_{x_{j,k} \in V(\mathcal{L})} x_{j,k}) \ y^{|E(\mathcal{L})|}\\
 \\
=& \sum_{\substack{ V \subset V(\mathcal{H}_{r,i}^{\infty}) \\ \text{Finite subset}}} \Big(\sum_{\substack{\mathcal{L} \in Sub_e(\mathcal{H}_{r,i}^{\infty})\\ V(\mathcal{L})=V}} (\prod_{x_{j,k} \in V}  x_{j,k}) \ y^{|E(\mathcal{L})|}\Big)\\\
  \\
  =&  \sum_{\substack{ V \subset V(\mathcal{H}_{r,i}^{\infty}) \\ \text{Finite subset}}} (\prod_{x_{j,k} \in V}  x_{j,k})  \Big(\sum_{\substack{\mathcal{L} \in Sub_e(\mathcal{H}_{r,i}^{\infty}) \\ V(\mathcal{L})=V}} y^{|E(\mathcal{L})|}\Big).\\
\end{align*}
On the one hand, as we mentioned before, for each finite sub-hypergraph $\mathcal{L'} \in Sub_v(\mathcal{H}_{r,i}^{\infty})$ with a finite $V(\mathcal{L'} )=V\subset  V(\mathcal{H}_{r,i}^{\infty})$, there exists a unique partition $\lambda\in \mathcal{N}_{r,i}$ such that $\mathcal{L'} =\mathcal{H}_{\lambda}$. Thus any sub-hypergraph $\mathcal{L} \in Sub_e(\mathcal{H}_{r,i}^{\infty})$ with $V(\mathcal{L} )=V$ is actually an induced edge sub-hypergraph of $\mathcal{L'}=\mathcal{H}_{\lambda}$ with $V(\mathcal{L} )=V.$ So we have:

$$S(\textbf{v},y)=\sum_{\lambda \in \mathcal{N}_{r,i}} (\prod_{x_{j,k} \in V(\mathcal{H}_{\lambda})}  x_{j,k}) \Big(\sum_{\substack{\mathcal{L} \in Sub_e(\mathcal{H}_{\lambda}) \\ V(\mathcal{L})=V(\mathcal{H_{\lambda}})}} (y)^{|E(\mathcal{L})|}\Big).$$ 

On the other hand, remember that the parts of $\lambda$ are in bijection with the vertices of $\mathcal{L}=\mathcal{H}_{\lambda}.$  By definition, $x_{j,k}$ is a vertices of $\mathcal{H}_{\lambda}$ if and only if $j$ is a part of $\lambda$ repeated at least $k$ times. This means that if we replace each $x_{j,k}$ by $q^j$ then 
 $\prod_{x_{j,k} \in V} x_{j,k} = q^{|\lambda|}$ and therefore:

 $$S_{r,i}^{w}(q,-1)=\sum_{\lambda \in \mathcal{N}_{r,i}} q^{|\lambda|} \Big(\sum_{\substack{\mathcal{L} \in Sub_e(\mathcal{H}_{\lambda}) \\ V(\mathcal{L})=V(\mathcal{H_{\lambda}})}} (-1)^{|E(\mathcal{L})|}\Big)=\sum_{\lambda \in \mathcal{N}_{r,i}} \delta(\lambda) q^{|\lambda|}.$$
\end{proof}

\section{Multigraded Hilbert series of an edge ideal of a hypergraph}
Let $\mathcal{H}$ be a simple hypergraph with a countable vertex set $V(\mathcal{H})=\{v_h| \ h \in I\}.$ In this section, with such $\mathcal{H},$ we associate a monomial ideal $\mathcal{I}(\mathcal{H})$ 
in a weighted polynomial ring $A$ whose variables are in bijection with the vertices of $\mathcal{H};$ we then 
consider a kind of Hilbert series $H_{\mathcal{H}}$ of this ideal, this is the generating series of the monomials 
in the quotient of $A$ by $ \mathcal{I}(\mathcal{H});$ this latter series is expressed in  Proposition \ref{Hilbert series} via $S(\textbf{v},y).$\\

Let $\mathbb{K}$ be a field of characteristic zero. Consider the polynomial ring $A=\mathbb{K}[v_h| \ h\in I].$ To each edge $e=(v_{h_1},\cdots,v_{h_\ell})$ of $\mathcal{H}$ we can associate a monomial $m_e=v_{h_1}\cdots v_{h_\ell}\in A.$
\\
\\
\textit{The edge ideal} $\mathcal{I}(\mathcal{H})$  of $\mathcal{H}$ is a square-free monomial ideal in $A$ whose generators are obtained from the edges of $\mathcal{H}.$ i.e.,

$$\mathcal{I}(\mathcal{H})=\langle m_{e_h}| \ e_h\in E(\mathcal{H}) \rangle.$$

\begin{defi} Let $\mathcal{H}$ be a simple hypergraph with a countable vertex set $V(\mathcal{H})=\{v_h| \ h \in I\}$. \textit{The Hilbert series} $H_{\mathcal{H}}=H_{A/ \mathcal{I} (\mathcal{H})}$ of $\mathcal{H}$ is the following multivariable series in $\textbf{v}:=(v_h)_{h \in I}.$:

$$H_{\mathcal{H}}(\textbf{v})=H_{A/ \mathcal{I} (\mathcal{H})}(\textbf{v})=\sum_{\substack{m\in A \setminus \mathcal{I}(\mathcal{H})  \\ \text{monomial}}} m$$
which is a series whose variables are the vertices of $\mathcal{H}.$
\end{defi}

\begin{prop}\label{Hilbert series} Let $\mathcal{H}$ be a simple hypergraph with a countable vertex set $V(\mathcal{H})=\{v_h| \ h \in I\}$. Then:

$$H_{\mathcal{H}}(\textbf{v})= \frac{S(\textbf{v},-1)}{\prod_{v_h \in V}(1-v_h)}.$$

\end{prop}

\begin{proof} Consider the hypergraph $\mathcal{L}$ with $V(\mathcal{L})=V(\mathcal{H})=V$ and $E(\mathcal{L})=\emptyset.$  We have:

$$H_{\mathcal{L}}(\textbf{v})=H_{A}(\textbf{v})=\sum_{\substack{m\in A  \\ \text{monomial}}} m = \sum_{i_h \in \mathbb{N}} (\prod_{v_h \in V} v_{h}^{i_h})
 = \prod_{v_h \in V} (\sum_{i_h \in \mathbb{N}} v_{h}^{i_h})
= \frac{1}{\prod_{v_h \in V}(1-v_h)}.$$

In order to compute $H_{\mathcal{H}}(\textbf{v})$ we have to consider all monomials of $A$ which are not in the monomial ideal  $\mathcal{I}(\mathcal{H})=\langle m_{e_h}| \ e_h\in E(\mathcal{H}) \rangle$.
Note that a monomial $m'\in A$ is in the ideal $\mathcal{I}(\mathcal{H})$ if and only if it is a multiple of \textbf{at least} one generator of $\mathcal{I}(\mathcal{H}).$ i.e.,
$$m'\in \mathcal{I}(\mathcal{H}) \iff \exists e_h \in E({\mathcal{H}}), \ \exists m \in A, \ m'=m_{e_h}.m.$$
So there could exist \textbf{several different generators} $m_{e_1},\cdots, m_{e_h}$ of $\mathcal{I}(\mathcal{H})$ such that $m'$ is a multiple of each of these generators.
 Thus $m'$ is also a multiple of the least common multiple of any subset of ${m_{e_1},\cdots, m_{e_h}}$.

 Therefore, if we denote by lcm$(m_{e_1},\cdots, m_{e_h})$ the least common multiple of 
 
 $m_{e_1},\cdots, m_{e_h}$ then, in order to remove all monomials of $\mathcal{I}(\mathcal{H})$ from $A$ \textbf{only once} we have to:
\begin{itemize}
\item[-] Add all monomials of $A$ (which is equivalent to compute $H_{A}(\textbf{v})=H_{\mathcal{L}}(\textbf{v})$).

\item[-] remove once the monomials of the from $m_e.m$ for some $e\in E(\mathcal{H})$ and $m\in A$. i.e,
$$- \sum_{\substack{m \in A \\ e \in E(\mathcal{H})}} m_e.m=- H_{A}(\textbf{v}) \sum_{e \in E(\mathcal{H})} m_e$$

\item[-] Since in the previous step we have removed twice the monomials of the form lcm$(m_{e_1},m_{e_2}).m$ so we have to add them ones. which means adding the following series:

$$ \sum_{\substack{m \in A \\ \{e_1,e_2\} \subset E(\mathcal{H})}}  \text{lcm}(m_{e_1},m_{e_2}).m= H_{A}(\textbf{v}) \sum_{ \{e_1,e_2\} \subset E(\mathcal{H})} \text{lcm}(m_{e_1},m_{e_2}).$$
\item[-] Once again, since in the previous step we have added the monomials of the form lcm$(m_{e_1},m_{e_2},m_{e_3}).m$  twice, we now need to remove them once and so on.  
\end{itemize}
Thus:
\begin{align*}
 H_{\mathcal{H}}(\textbf{v})=H_{A}(\textbf{v})&  \Big(1- \sum_{e \in E(\mathcal{H})} m_e+ \sum_{ \{e_1,e_2\} \subset E(\mathcal{H})} \text{lcm}(m_{e_1},m_{e_2})\\\
 \\
  & -\sum_{ \{e_1,e_2,e_3\} \subset E(\mathcal{H})} \text{lcm}(m_{e_1},m_{e_2},m_{e_3})+\cdots+   .\\
  &+(-1)^k\sum_{ \{e_1,\cdots,e_k\} \subset E(\mathcal{H})} \text{lcm}(m_{e_1},\cdots,m_{e_k})+\cdots \Big) \\
\end{align*}
If we denote by $T_k=(-1)^k\sum_{ \{e_1,\cdots,e_k\} \subset E(\mathcal{H})} \text{lcm}(m_{e_1},\cdots,m_{e_k})$ then we have:

$$H_{\mathcal{H}}(\textbf{v})=H_{A}(\textbf{v})(1+T_1+\cdots+T_k+\cdots).$$

Note that choosing $k$ edges $\{e_1,\cdots,e_k\} \subset E(\mathcal{H})$ is equivalent to considering an edge induced sub-hypergraph $\mathcal{L}$ of $\mathcal{H}$ with $|E(\mathcal{L})|=k$ whose edge set is  $E(\mathcal{L})=\{e_1,\cdots,e_k\} \subset E(\mathcal{H})$ and whose vertex set is the union of the vertices of its edges. Thus lcm$(m_{e_1},\cdots,m_{e_k})$ is equal to the product of the vertices of $\mathcal{L}$ and therefore:
$$T_k= \sum_{\substack{\mathcal{L} \subset Sub_e(\mathcal{H}) \\ |E(\mathcal{L})|=k}} (-1)^{|E(\mathcal{L})|} (\prod_{v_{\ell} \in V(\mathcal{L})}v_{\ell}),$$
and 
\begin{align*}
H_{\mathcal{H}}(\textbf{v})=H_{A}(\textbf{v})&\Big(1+\sum_{k \in \mathbb{N}^*}\sum_{\substack{\mathcal{L} \subset Sub_e(\mathcal{H}) \\ |E(\mathcal{L})|=k}} (-1)^{|E(\mathcal{L})|} (\prod_{v_{\ell} \in V(\mathcal{L})}v_{\ell})\Big)\\
 =H_{A}(\textbf{v}) &\Big( 1+\sum_{\substack{\mathcal{L} \subset Sub_e(\mathcal{H}) \\ \text{finite}}} (-1)^{|E(\mathcal{L})|}  (\prod_{v_{\ell} \in V(\mathcal{L})}v_{\ell})\Big)\\
  =H_{A}(\textbf{v}) & \ S(\textbf{v},-1)\\
  \\
  =& \frac{S(\textbf{v},-1)}{\prod_{v_h \in V}(1-v_h)}.
\end{align*}

\end{proof}

\section{Hilbert-Poincar\'e series and signature of the $(r,i)$-Neighborly partition}

For the integers $1\leq i \leq r$ consider the $\mathbb{K}$-algebra: 
$$\mathcal{P}_{r,i}=\mathbb{K}[x_{j,k}| \ x_{j,k}\in V(\mathcal{H}_{r,i}^{\infty})]/\mathcal{I}(\mathcal{H}^{\infty}_{r,i}),$$
which is graded by giving the weight $j$ to $x_{j,k}.$ This means that there exist finite-dimensional algebras $\mathcal{P}_{(r,i),n},$ each generated by the monomials of weight $n$ which do not belong to the ideal $\mathcal{I}(\mathcal{H}^{\infty}_{r,i})$ and such that  $\mathcal{P}_{r,i}=\oplus_{n \in \mathbb{N}} \mathcal{P}_{(r,i),n}.$ By definition, the \textit{Hilbert-poincar\'e series} of $\mathcal{P}_{r,i}$ is given by:

$$HP_{\mathcal{P}_{r,i}}(q)=\sum_{n\in \mathbb{N}} \dim_{\mathbb{K}}(\mathcal{P}_{(r,i),n})q^n.$$

Using Proposition \ref{Hilbert series} and Lemma \ref{S-series} we can see the relation between this $q$ series and the signature of the neighborly partitions in the following proposition:

\begin{prop}\label{HP} For the integers $1\leq i \leq r$ we have:
$$HP_{\mathcal{P}_{r,i}}(q)=\frac{\sum_{\lambda \in \mathcal{N}_{r,i}} \delta(\lambda) q^{|\lambda|}}{(1-q)^i \prod_{j \in \mathbb{N}\setminus \{0,1\}}(1-q^j)^r}.$$
\end{prop}
\begin{proof}

Note that the monomials of weight $n$ in $\mathbb{K}[x_{j,k}| \ x_{j,k}\in V(\mathcal{H}_{r,i}^{\infty})]$ are those whose sum of the first indices is equal to $n.$ i.e., a monomial $m=x_{j_1,k_1}\cdots x_{j_{\ell} k_\ell}$ is of weight $n$ if $\sum_{s=1}^{\ell} j_s=n.$ This means that if we replace $x_{j_s,k_s}$ by $q^{j_s}$ in $m$ then we obtain $q^n.$ Thus, by definition of the Hilbert series of a hypergraph and the Hilbert-Poincar\'e series of a graded Algebra, if we denote by $H_{\mathcal{H}_{r,i}}^w(q)$ the series obtained from the Hilbert series of $\mathcal{H}_{r,i}^{\infty}$ then we have:

$$HP_{\mathcal{P}_{r,i}}(q)=H_{\mathcal{H}_{r,i}}^w(q).$$

Applying Proposition \ref{Hilbert series} and Lemma \ref{S-series} on $H_{\mathcal{H}_{r,i}}^w(q)$ we obtain:

$$HP_{\mathcal{P}_{r,i}}(q)= \frac{S_{r,i}^{w}(q,-1)}{(1-q)^i \prod_{j \in \mathbb{N}\setminus \{0,1\}}(1-q^j)^r}=\frac{\sum_{\lambda \in \mathcal{N}_{r,i}} \delta(\lambda) q^{|\lambda|}}{(1-q)^i \prod_{j \in \mathbb{N}\setminus \{0,1\}}(1-q^j)^r}.$$

\end{proof}

\section{Gordon's identities and $(r,i)$-Neighborly partitions}
In this section we give the relation between the partitions appearing in the Gordon's identities and the $(r,i)$-Neighborly partitions. Using this relation, we  prove our main theorem. 
\\
To do so, we need to use the \textit{polarization} of the monomial ideals. This is an operation which allows us to obtain a square-free monmial ideal from a non-square-free one by turning the monomial $x_\ell^s \in \mathbb{K}[x_j|\ j \in B]$ to the square free monomial $x_{\ell,1}\cdots x_{\ell,s} \in \mathbb{K}[x_{j,k}|\ j\in B, k \in C].$
There is a close relation between the Hilbert-Poincar\'e series of the quotient ring of an ideal and the Hilbert-Poincar\'e series of the quotient ring of its polarization ideal (See Corollary $1.6.3$ of \cite{HH}, see also \cite{P} and \cite{MS}).

 Denote by $J_{r,i}$ the following monomial ideal:
$$J_{r,i}=\langle x_{1}^i, x_{\ell}^{s}x_{\ell+1}^{r-s} | \ 1 \leq s \leq r\ \text{and} \ \ell \geq 1 \rangle \subset \mathbb{K}[x_j| \ j\geq 1].$$

Note that the $\mathcal{I}(\mathcal{H}^{\infty}_{r,i}),$ is the polarization of the ideal $J_{r,i}$. In \cite{BMS} (see also \cite{BMS1}) Bruschek, Mourtada and Schepers proved that if we give to $x_j$ the weight $j$ then Hilbert-poincar\'e series of the graded algebra $\mathbb{K}[x_j | \ j\geq 1]/ J_{r,i}$ is  equal to the generating series of partitions appearing in the Gordon's identities, which are counted by $B_{r,i}(n).$ We use this result and the relation between $\mathcal{I}(\mathcal{H}^{\infty}_{r,i})$ and $J_{r,i}$ to prove our main results:

\begin{theo}\label{generatrice} Let $1\leq i \leq r$ be the integers. We have 

$$\sum_{\lambda \in \mathcal{N}_{r,i}}\delta(\lambda) q^{|\lambda|}= \sum_{n\in \mathbb{N}}\mathcal{R}_{r,i}(n)q^n= \prod_{j \equiv 0, \pm i [2r+1]} (1-q^j),$$
where $|\lambda|$ is the sum of the parts of $\lambda$.
\end{theo}
\begin{proof}
On the one hand, since $\mathcal{I}(\mathcal{H}^{\infty}_{r,i})$ is the polarization of the ideal $J_{r,i}$, we have the following relation between the Hilbert-Poincar\'e series of the quotient rings of these ideals (See Corollary $1.6.3$ of \cite{HH}, see also \cite{P} and \cite{MS}):
$$
HP_{\mathcal{P}_{r,i}}(q)= \frac{HP_{ \frac{\mathbb{K}[x_j| \ j\geq 1]}{J_{r,i}}}(q)}{(1-q)^{i-1} \prod_{j \in \mathbb{N}\setminus \{0,1\}}(1-q^j)^{r-1}}
$$
On the other hand, since the Hilbert-poincar\'e series of $\mathbb{K}[x_j | \ j\geq 1]/ J_{r,i}$  is equal to the generating series of partitions appearing in the Gordon's identities, we have:
 $$HP_{ \frac{\mathbb{K}[x_j| \ j\geq 1]}{J_{r,i}}}(q)=\sum_{n\in \mathbb{N}} B_{r,i}(n)q^n=\sum_{n\in \mathbb{N}} A_{r,i}(n)q^n=\frac{ \prod_{j \equiv 0, \pm i [2r+1]} (1-q^j)}{\prod_{j\in \mathbb{N}}(1-q^j)}.$$
  and thus:
 $$
HP_{\mathcal{P}_{r,i}}(q)= \frac{ \prod_{j \equiv 0, \pm i [2r+1]} (1-q^j)}{(1-q)^{i} \prod_{j \in \mathbb{N}\setminus \{0,1\}}(1-q^j)^{r}}.
$$
 
 By comparing this formula with the formula in Proposition \ref{HP} we obtain:
 
$$
\sum_{\lambda \in \mathcal{N}_{r,i}} \delta(\lambda) q^{|\lambda|}= \prod_{j \equiv 0, \pm i [2r+1]} (1-q^j).
$$
The right-hand side of the equation above is the generating series of $\mathcal{R}_{r,i}(n)$.
\end{proof}

As we mentioned in the introduction, Theorem \ref{generatrice} is equivalent to the following theorem:

\begin{theo}\label{coeff.} Let $1\leq i \leq r$ be the integers. Then:
$$\sum_{\lambda \in \mathcal{N}_{r,i}(n)} \delta(\lambda)=\sum_{\lambda\in \mathcal{R}_{r,i}(n)} (-1)^{\ell(\lambda)}.$$
\end{theo}

\begin{proof}
Fix an integer $n \in \mathbb{N}$. In order to prove this theorem we prove that the coefficients of $q^n$ in both sides of the equation in Theorem \ref{generatrice} are equal.
\\
On the one hand, we have:
 $$\sum_{\lambda \in \mathcal{N}_{r,i}} \delta(\lambda) q^{|\lambda|}=\sum_{n\in \mathbb{N}} (\sum_{\lambda \in \mathcal{N}_{r,i}(n)} \delta(\lambda)) q^n.$$
 
On the other hand, in order to find the coefficient of $q^n$ in $ \prod_{j \equiv 0, \pm i [2r+1]} (1-q^j)$,  we take a partition $\lambda:(\lambda_1, \cdots, \lambda_{\ell})\in \mathcal{R}_{r,i}(n)$ of length $\ell$. Thus, for each $1\leq j \leq \ell$ we know that $\lambda_j \equiv 0,\pm i [2r+1]$ with $\sum_{j=1}^{\ell}\lambda_j=n$ and we have:
 $$1.(-q^{{\lambda_1}})(-q^{{\lambda_2}})\cdots(-q^{{\lambda_{\ell}}})=(-1)^{\ell}q^{\sum_{j=1}^{\ell}{\lambda_j}}=(-1)^{\ell}q^{n}.$$
This proves that the coefficient of $q^n$ in the right-hand side of the equation in Theorem \ref{generatrice} is equal to:
$$\sum_{\lambda\in \mathcal{R}_{r,i}(n)} (-1)^{\ell(\lambda)}.$$
Thus:
$$\sum_{\lambda \in \mathcal{N}_{r,i}(n)} \delta(\lambda)=\sum_{\lambda\in \mathcal{R}_{r,i}(n)} (-1)^{\ell(\lambda)}.$$
\end{proof}





\bibliographystyle{acm}
\bibliography{Afsharijoo_Mourtada_Neighborly}

\hspace{1cm}

\address{Instituto de Matem\'atica Interdisciplinar, Departamento de \'Algebra, Geometr\'ia y Topolog\'ia,
Facultad de Ciencias Matem\'aticas, Universidad Complutense de Madrid, Plaza de las Ciencias 3,
Madrid 28040, Espa\~na.\\
\email{pafshari@ucm.es} \\}

\address{Universit\'e Paris Cit\'e, Sorbonne Universit\'e, CNRS, Institut de Math\'ematiques de Jussieu-Paris Rive Gauche, Paris, 75013, France,\\
\email{hussein.mourtada@imj-prg.fr} }

\end{document}